\documentclass[11pt,a4,fleqn]{article}
\usepackage{graphicx}
\usepackage{amsmath,amssymb,latexsym,graphics,epsfig}
\usepackage{hyperref}
\usepackage{color}
\usepackage{amsthm}

\newtheorem{theorem}{\bf Theorem}[section]
\newtheorem{proposition}[theorem]{\bf Proposition}
\newtheorem{lemma}[theorem]{\bf Lemma}
\newtheorem{corollary}[theorem]{\bf Corollary}

\newcommand{\lal}{\lambda_{\ell}}
\newcommand{\mul}{\mu_{\ell}}
\newcommand{\nul}{\nu_{\ell}}
\newcommand{\G}{\Gamma}

\newcommand{\pf}{\noindent{\em Proof: }}
\newcommand{\epf}{\hfill\hbox{\rule{3pt}{6pt}}\\}

\numberwithin{equation}{section}

\begin{document}
\title{{\Large Strongly walk-regular graphs}\footnote{This version is published in Journal of Combinatorial Theory, Series A
120 (2013), 803--810.}}

\author{E.R. van Dam$^{\textrm{a}}$   \quad  G.R. Omidi$^{\textrm{b},\textrm{c},1}$  \\[2pt]
{\small  $^{\textrm{a}}$Tilburg University, Dept. Econometrics and Operations Research,}\\
{\small  P.O. Box 90153, 5000\,LE, Tilburg, The Netherlands}\\
{\small  $^{\textrm{b}}$Dept. Mathematical Sciences, Isfahan University of Technology},\\
{\small Isfahan, 84156-83111, Iran}\\
{\small $^{\textrm{c}}$School of Mathematics, Institute for Research in Fundamental Sciences (IPM),}\\
{\small P.O. Box 19395-5746, Tehran, Iran }\\[2pt]
{edwin.vandam@uvt.nl \quad romidi@cc.iut.ac.ir}}

\date{}

\maketitle \footnotetext[1] {This research is partially
carried out in the IPM-Isfahan Branch and in part supported
by a grant from IPM (No. 91050416).} \vspace*{-0.5cm}

\begin{abstract}
\noindent We study a generalization of strongly regular graphs. We call a graph
strongly walk-regular if there is an $\ell >1$ such that the number of walks of
length $\ell$ from a vertex to another vertex depends only on whether the two
vertices are the same, adjacent, or not adjacent. We will show that a strongly
walk-regular graph must be an empty graph, a complete graph, a strongly regular
graph, a disjoint union of complete bipartite graphs of the same size and
isolated vertices, or a regular graph with four eigenvalues. Graphs from the
first three families in this list are indeed strongly $\ell$-walk-regular for
all $\ell$, whereas the graphs from the fourth family are $\ell$-walk-regular for every odd $\ell$.
The case of regular graphs with four eigenvalues is the most interesting (and
complicated) one. Such graphs cannot be strongly $\ell$-walk-regular for even
$\ell$. We will characterize the case that regular four-eigenvalue graphs are
strongly $\ell$-walk-regular for every odd $\ell$, in terms of the eigenvalues.
There are several examples of infinite families of such graphs. We will show
that every other regular four-eigenvalue graph can be strongly
$\ell$-walk-regular for at most one $\ell$. There are several examples of
infinite families of such graphs that are strongly $3$-walk-regular. It however
remains open whether there are any graphs that are strongly $\ell$-walk-regular
for only one
particular $\ell$ different from $3$.\\


\end{abstract}

\section{Introduction}

A strongly regular graph is a regular graph such that the number of common neighbors of two distinct vertices depends
only on whether these vertices are adjacent or not. Strongly regular graphs are well-studied combinatorial and
algebraic objects that arise in the study of finite groups, codes, designs, and finite geometries (see the book by
Brouwer and Haemers \cite{bh12}, for example). Strongly regular graphs also have prominent generalizations in the form
of distance-regular graphs (see the book by Brouwer, Cohen, and Neumaier \cite{bcn89}, or the recent survey
\cite{dkt12}) and association schemes (see the recent survey by Martin and Tanaka \cite{martintanaka}, for example).

Because the number of common neighbors of two vertices in a graph is the same
as the number of walks of length two from one vertex to the other, we can
generalize the concept of strongly regular graph by considering the number of
walks of length --- say --- $\ell$ from one vertex to the other. Thus, we call
a graph {\em strongly walk-regular} if there is an $\ell >1$ such that the
number of walks of length $\ell$ from a vertex to another vertex depends only
on whether the two vertices are the same, adjacent, or not adjacent.

We remark that Godsil and McKay \cite{gmk} introduced walk-regular graphs as those graphs for which the number of walks
of length $\ell$ from a vertex to itself does not depend on the chosen vertex, but only on $\ell$. Rowlinson \cite{r97}
characterized distance-regular graphs by the property that the number of walks of length $\ell$ between two vertices
depends only on $\ell$ and the distance between the two vertices. Dalf\'{o}, Van Dam, Fiol, Garriga, and Gorissen
\cite{ddfgg10} studied several versions of so-called almost distance-regularity by generalizing Rowlinson's
characterization. Let us emphasize that in these three papers \cite{ddfgg10, gmk, r97}, the mentioned properties should
(typically) hold for {\em all} $\ell$, whereas here (in this paper) we require a certain property to hold for only one
specific $\ell$.

At first sight, strongly walk-regular graphs also seem to be related to
Kotzig's conjecture about generalized friendship graphs. The friendship theorem
states that if a graph has the property that every pair of vertices are
connected by a unique {\em path} of length $2$, then there is a vertex that is
adjacent to every other vertex; the proof of this result is in ``THE BOOK''
\cite[Ch. 34]{BOOK}. Kotzig \cite{Kotzig} conjectured that for $\ell>2$, there
are no graphs with the property that every pair of (distinct) vertices are
connected by a unique path of length $\ell$. According to Kostochka (see
\cite[p. 225]{BOOK}) this has been proven for $\ell \leq 33$. This conjecture
and our problem are however very different. Not only do they differ because of
the difference between paths and walks, but also in the fact that Kotzig only
considers pairs of distinct vertices. We therefore leave to the reader the study of graphs with
the property that the number of {\em paths} of length $\ell$ between two distinct vertices
depends only on whether the two vertices are adjacent or not.

In this paper, we will show that a strongly walk-regular graph must be an empty graph, a
complete graph, a strongly regular graph, a disjoint union of complete
bipartite graphs of the same size and isolated vertices, or a regular graph
with four eigenvalues. This conclusion can be drawn from Theorem \ref{thm:3or4} and Proposition \ref{prop:mu=0}.
Graphs from the first three families in this list are
indeed strongly $\ell$-walk-regular for all $\ell$, whereas the disjoint unions
of complete bipartite graphs of the same size and isolated vertices are
$\ell$-walk-regular for every odd $\ell$. The case of regular graphs with four
eigenvalues, which is discussed in Section \ref{sec:4ev}, is the most interesting
(and complicated) one. Such graphs cannot be strongly $\ell$-walk-regular for even $\ell$.
We will characterize the case
that regular four-eigenvalue graphs are strongly $\ell$-walk-regular for every
odd $\ell$, in terms of the eigenvalues. There are several examples of infinite
families of such graphs. We will show that every other regular four-eigenvalue
graph can be strongly $\ell$-walk-regular for at most one $\ell$. There are
several examples of infinite families of such graphs that are strongly
$3$-walk-regular. It however remains open whether there are any graphs that are
strongly $\ell$-walk-regular for only one particular $\ell$ different from $3$.

\section{The definition and some basic observations}

We will now give the formal definition of the concept of a strongly walk-regular graph as mentioned in the
introduction. For $\ell >1$, a graph $\G$ is a {\em strongly $\ell$-walk-regular graph with parameters}
$(\lal,\mul,\nul)$ if there are $\lal$, $\mul$, and $\nul$ walks of length $\ell$ between every two adjacent, every two
non-adjacent, and every two identical vertices, respectively. So indeed, every strongly regular graph with parameters
$(v,k,\lambda,\mu)$ is a strongly $2$-walk-regular graph with parameters $(\lambda,\mu,k)$. We exclude the case
$\ell=1$ because every graph would be a strongly $1$-walk-regular graph with parameters $(1,0,0)$. Besides strongly
regular graphs, also the empty and complete graphs give examples. Indeed, these graphs are clearly strongly
$\ell$-walk-regular for every $\ell$.

Let $A$ be the adjacency matrix of $\G$. A crucial observation is that $\G$ is
a strongly $\ell$-walk-regular graph if and only if $A^{\ell}$ is in the span
of $A,I$, and the all-ones matrix $J$. This follows from the following lemma
that translates the combinatorial definition into an algebraic formulation.
(For background on algebraic graph theory we refer to the book by Brouwer and
Haemers \cite{bh12}.)

\begin{lemma}\label{lem:A} Let $\ell>1$, and let $\G$ be a graph with adjacency matrix $A$.
Then $\G$ is a strongly $\ell$-walk-regular graph with parameters
$(\lal,\mul,\nul)$ if and only if $A^{\ell}+(\mul-\lal)A+(\mul-\nul)I=\mul J$.
\end{lemma}

\pf This follows by observing that $(A^{\ell})_{ij}$ counts the number of walks
of length $\ell$ between vertices $i$ and $j$. \epf

For a strongly regular graph $\G$, its adjacency algebra (that is, the algebra
spanned by all powers of $A$) equals $\langle A,I,J \rangle$, so this implies
that $\G$ is a strongly $\ell$-walk-regular graph for every $\ell>1$.

\begin{proposition}\label{prop:srg} Let $\G$ be a strongly regular graph.
Then $\G$ is a strongly $\ell$-walk-regular graph with parameters $(\lal,\mul,\nul)$
for every $\ell>1$ and some $\lal,\mul,$ and $\nul$.
\end{proposition}

At this point it seems appropriate to remark that Abreu, Funk, Labbate, and
Napolitano \cite{ghs} recently generalized the (strongly regular) Moore graphs
with girth $5$ by considering the equation $((k-1)I+J-A^2)^m=A$. This is
however a quite different equation in terms of $A,I,$ and $J$ than the one we
are considering.

Let us now make some more elementary observations that will help to see which
other graphs can be strongly walk-regular.

\begin{lemma} Let $\ell>1$, and let $\G$ be a
strongly $\ell$-walk-regular graph with parameters $(\lal,\mul,\nul)$ with
$\mul>0$. Then $\G$ is regular and connected.
\end{lemma}

\pf Because $\mul>0$, it follows from Lemma \ref{lem:A} that $J$ can be
expressed as a polynomial in $A$. Hoffman \cite{hoffman63} showed that this
implies that $\G$ is regular (by observing that $A$ commutes with $J$) and
connected. \epf

Note that also if $\mul=0$, the graph can be regular and connected. The
complete bipartite graph $K_{m,m}$ for example is strongly $3$-walk-regular
with parameters $(m^2,0,0)$. We will look further into the case $\mul=0$ in
Section \ref{sec:mu0}.

\section{Connected regular graphs}

We now first focus on the generic case and assume that the graphs are connected
and regular.

\begin{proposition}\label{eigenvalues}
A connected $k$-regular graph $\G$ on $v$ vertices is strongly
$\ell$-walk-regular with parameters $(\lal,\mul,\nul)$ if and only if all
eigenvalues except $k$ are roots of the equation
$$x^{\ell}+(\mul-\lal)x+\mul-\nul=0$$
and $k$ satisfies the equation
$$k^{\ell}+(\mul-\lal)k+\mul-\nul=\mul v.$$
\end{proposition}

\pf If $\G$ is a strongly $\ell$-walk-regular graph with parameters
$(\lal,\mul,\nul)$, then $A^{\ell}+(\mul-\lal)A+(\mul-\nul)I=\mul J$. Every
eigenvalue $\theta$ unequal to $k$ has an eigenvector that is orthogonal to the
all-ones vector (the eigenvector for eigenvalue $k$), hence by multiplying the
matrix equation with this eigenvector, we obtain the equation
$\theta^{\ell}+(\mul-\lal)\theta+\mul-\nul=0.$ By multiplying the matrix
equation with the all-ones vector, the equation
$k^{\ell}+(\mul-\lal)k+\mul-\nul=\mul v$ is derived.

On the other hand, consider the distinct eigenvalues $k=\theta_0>\theta_1> \cdots >\theta_d$ and the Hoffman polynomial
\cite{hoffman63} $h$ of $\G$. Here we normalize $h$ such that it is monic, so that
\begin{equation}\label{eq:hoffmanpol}
h(x)=\prod_{i=1}^d (x-\theta_i).
\end{equation}
In fact, $h(A)=\frac{h(k)}vJ$, because \eqref{eq:hoffmanpol} implies that $h(A)$ vanishes on each eigenspace of $A$
corresponding to an eigenvalue unequal to $k$, and thus it has rank 1 (because $\G$ is connected) with nontrivial
eigenvalue $h(k)$ having the all-ones vector as eigenvector.

Now, if the eigenvalue $\theta_i$ satisfies the equation $x^{\ell}+(\mul-\lal)x+\mul-\nul=0$ for every $i \neq 0$, then
the Hoffman polynomial divides $x^{\ell}+(\mul-\lal)x+\mul-\nul$. Because $h(A)=\frac{h(k)}vJ$ and $p(A)J=p(k)J$ for
every polynomial $p$, this implies that $A^{\ell}+(\mul-\lal)A+(\mul-\nul)I=tJ$ for some $t$. As before, this implies
that $k^{\ell}+(\mul-\lal)k+\mul-\nul=t v$. But if $k^{\ell}+(\mul-\lal)k+\mul-\nul=\mul v$, then $t=\mul$, and so $\G$
is a strongly $\ell$-walk-regular graph with parameters $(\lal,\mul,\nul)$. \epf

As a corollary, we obtain the following characterization that we will use later
on.

\begin{corollary}\label{cor:eigenvalues} A connected regular graph is strongly
$\ell$-walk-regular if and only if its Hoffman polynomial divides the
polynomial $x^{\ell}+ex+f$ for some $e$ and $f$.
\end{corollary}

\pf One direction follows immediately from Proposition \ref{eigenvalues}. To
show the other direction, assume that the Hoffman polynomial divides the
polynomial $x^{\ell}+ex+f$ for some $e$ and $f$. Let $v$ be the number of
vertices and $k$ the valency of the graph. Clearly, all eigenvalues besides $k$
are roots of $x^{\ell}+ex+f$. Now let $\mul=\frac1v(k^{\ell}+ek+f),
\lal=\mul-e,$ and $\nul=\mul-f$. Then Proposition \ref{eigenvalues} implies
that the graph is $\ell$-walk-regular with parameters $(\lal,\mul,\nul)$. \epf

By bounding the number of real roots of $x^{\ell}+ex+f$, we can bound the
number of distinct eigenvalues of a strongly walk-regular graph.

\begin{lemma}\label{lem:descartes} Let $\ell>1$, and let $p(x)=x^{\ell}+ex+f$ for some $e$ and
$f$. Then $p$ has at most three real roots. If $\ell$ is even, then $p$ has at
most two real roots.
\end{lemma}

\pf The derivative of $p$ equals $p'(x)=\ell x^{\ell-1}+e$. So $p'$ has at most
two real roots if $\ell$ is odd, and at most one real root if $\ell$ is even.
From this, the result follows. \epf

We remark that if we know the sign of $e$ ($=\mul-\lal$) and $f$
($=\mul-\nul$), then some more information about the sign and number of real
roots can be derived by using Descartes' rule of signs (see \cite{descartes}).

Because eigenvalues of graphs are real, the following result follows
immediately from Corollary \ref{cor:eigenvalues} and Lemma \ref{lem:descartes}.

\begin{theorem}\label{thm:3or4} Let $\G$ be a connected regular
strongly $\ell$-walk-regular graph with $\ell>1$. Then $\G$ has at most four
distinct eigenvalues. Moreover, if $\ell$ is even, then $\G$ has at most three
distinct eigenvalues and hence $\G$ is either a complete graph or a strongly
regular graph.
\end{theorem}

\section{Graphs with four eigenvalues}\label{sec:4ev}

Let us continue with the case of connected regular graphs. By Proposition \ref{prop:srg} and Theorem \ref{thm:3or4}, it
remains to determine the strongly $\ell$-walk-regular graphs with four (distinct) eigenvalues for odd $\ell$. Regular
four-eigenvalue graphs have been studied by Van Dam and Spence \cite{vD-4e,vDS}. It is for example known that such
graphs are walk-regular in the sense of Godsil and McKay \cite{gmk}. Because strongly regular graphs are also
walk-regular, it follows that every connected regular strongly walk-regular graph is indeed walk-regular. We remark
that in Section \ref{sec:mu0}, we will see some {\em nonregular} strongly walk-regular graphs, and these graphs cannot
be walk-regular, because every walk-regular graph is regular.

Let $\G$ be a connected $k$-regular graph with four distinct eigenvalues
$k>\theta_1>\theta_2>\theta_3$. Then by working out the Hoffman polynomial
\eqref{eq:hoffmanpol}, it follows that
\begin{equation}\label{eq:Hoffman}
A^3-(\theta_1+\theta_2+\theta_3)A^2+(\theta_1\theta_2+\theta_1\theta_3+\theta_2\theta_3)A-\theta_1\theta_2\theta_3 I=t J
\end{equation}
for some $t$ (this $t$ can be expressed in terms of the number of vertices and
the eigenvalues, but this doesn't play any role in the following). This proves
the following.

\begin{proposition}\label{prop:3swr} Let $\G$ be a connected regular graph with
four distinct eigenvalues $k>\theta_1>\theta_2>\theta_3$. Then $\G$ is strongly
$3$-walk-regular if and only if $\theta_1+\theta_2+\theta_3=0$.
\end{proposition}

\pf This follows from \eqref{eq:Hoffman}, and the fact that $A^2$ cannot be
expressed in terms of $A,I,J$, because otherwise $\G$ would be strongly
regular, and thus have three distinct eigenvalues. \epf

Several infinite families of examples of regular four-eigenvalue graphs with $\theta_1+\theta_2+\theta_3=0$ can be
constructed using the methods of Van Dam \cite{vD-4e}. Many, but certainly not all of these examples have integer
eigenvalues only. For an example with noninteger eigenvalues, consider a conference graph, that is, a strongly regular
graph with parameters $(4\mu+1,2\mu,\mu-1,\mu)$. Its 3-clique extension (that is, every vertex is replaced by a clique
of size 3, and vertices of different cliques are adjacent if and only if the original vertices are adjacent) has
distinct eigenvalues $k, \frac{3}{2}(-1+\sqrt{4\mu+1})+2,-1$, and $\frac{3}{2}(-1-\sqrt{4\mu+1})+2$.

There are also some distance-regular graphs with diameter three that are
strongly $3$-walk-regular, such as the line graph of the incidence graph of the
Fano plane, the Perkel graph, and the Hamming graph $H(3,3)$.

Another family of examples is given by the complements of the graphs $K_{m,m}\oplus K_m$ (the sum \cite{cds82} --- also
called Cartesian product --- of a complete bipartite graph and a complete graph, see \cite{vD-4e}). The distinct
eigenvalues of the complement of $K_{m,m}\oplus K_m$ are $k=2m(m-1),m,0,$ and $-m$. The zero eigenvalue simplifies
\eqref{eq:Hoffman} even further, and it implies that the graph is strongly $\ell$-walk-regular for every odd $\ell$.

\begin{proposition} Let $\G$ be a connected regular graph with four distinct
eigenvalues $k>\theta_1>\theta_2>\theta_3$. If $\theta_2=0$ and
$\theta_3=-\theta_1$, then $\G$ is strongly $\ell$-walk-regular for every odd
$\ell$.
\end{proposition}

\pf If $\theta_2=0$ and $\theta_3=-\theta_1$, then the three eigenvalues
$\theta_1,\theta_2,$ and $\theta_3$ are roots of the polynomial
$x^{\ell}-\theta_1^{\ell-1}x$ for every odd $\ell$, so $\G$ is
$\ell$-walk-regular for every odd $\ell$ by Corollary \ref{cor:eigenvalues}.
\epf

On the other hand, if $\theta_2 \neq 0$ or $\theta_3 \neq -\theta_1$, then $\G$
is $\ell$-walk-regular for at most one $\ell>1$. In order to show this, we will
use the following characterization.

\begin{proposition} Let $\G$ be a connected regular graph with four distinct
eigenvalues $k>\theta_1>\theta_2>\theta_3$ and let $\ell \geq 3$. Then $\G$ is
strongly $\ell$-walk-regular if and only if
\begin{equation}\label{eq:lswr}
(\theta_2-\theta_3)\theta_1^{\ell}+(\theta_3-\theta_1)\theta_2^{\ell}+(\theta_1-\theta_2)\theta_3^{\ell}=0.
\end{equation}
\end{proposition}

\pf As before, we will use that $\G$ is strongly $\ell$-walk-regular if and only if the Hoffman polynomial
$x^3+ax^2+bx+c$ divides a polynomial of the form $x^{\ell}+ex+f$. Let $x^{\ell-3}+\sum_{i=3}^{\ell-1}\alpha_i
x^{\ell-1-i}$ be a putative quotient of these two polynomials (the unusual indexing will turn out to be convenient),
and also define $\alpha_0=\alpha_1=0$ and $\alpha_2=1$ (these can also be interpreted as coefficients of the same
polynomial). Now consider the product of this polynomial and the Hoffman polynomial. The coefficients of $x^{\ell-1}$
down to $x^3$ in this product vanish if and only if the coefficients $\alpha_i$ satisfy the recurrence relation
$\alpha_i+a\alpha_{i-1}+b\alpha_{i-2}+c\alpha_{i-3}=0$ for $i=3,4,\dots,\ell-1$. The coefficient of $x^2$ in the
product equals $a\alpha_{\ell-1}+b\alpha_{\ell-2}+c\alpha_{\ell-3}$, so if we extend the recurrence by one more step by
defining $\alpha_{\ell}$ such that $\alpha_{\ell}+a\alpha_{\ell-1}+b\alpha_{\ell-2}+c\alpha_{\ell-3}=0$, then it
follows that $\G$ is strongly $\ell$-walk-regular if and only if $\alpha_{\ell}=0$.

Now all that is left is to find the solution for the recurrence relation. But
the characteristic polynomial of this recurrence relation is the Hoffman
polynomial that has distinct roots $\theta_1,\theta_2$, and $\theta_3$.
Therefore $\alpha_i=\beta_1 \theta_1^i+\beta_2 \theta_2^i+\beta_3 \theta_3^i$
for $i=0,1,\dots,\ell$ for some $\beta_1,\beta_2,$ and $\beta_3$. Given that
$\alpha_0=\alpha_1=0$, and $\alpha_2=1$, it follows that
$$-(\theta_1-\theta_2)(\theta_2-\theta_3)(\theta_3-\theta_1)\alpha_i=(\theta_2-\theta_3)\theta_1^{i}+(\theta_3-\theta_1)\theta_2^{i}+(\theta_1-\theta_2)\theta_3^{i},$$
which finishes the proof. \epf

\begin{theorem}
Let $\G$ be a connected regular graph with four distinct eigenvalues
$k>\theta_1>\theta_2>\theta_3$. If $\theta_2 \neq 0$ or $\theta_3 \neq
-\theta_1$, then there is at most one $\ell>1$ such that $\G$ is strongly
$\ell$-walk-regular.
\end{theorem}

\pf Suppose on the contrary that $\G$ is strongly $\ell$-walk-regular and
strongly $m$-walk-regular with $m>\ell>1$. Note that \eqref{eq:lswr} implies
that $\theta_2 \neq 0$ if and only if $\theta_3 \neq -\theta_1$, so we may
assume both. Similarly, we may assume that $\theta_3 \neq -\theta_2$. A final
assumption that we will make is that $\theta_1>|\theta_2|$ and
$\theta_1>|\theta_3|$. We can do this without loss of generality because if we
replace the eigenvalues $\theta_1,\theta_2,$ and $\theta_3$ by their opposites,
then $\eqref{eq:lswr}$ still holds (and we will not make use of other specific
properties of the eigenvalues).

Note that \eqref{eq:lswr} trivially holds for $\ell=1$, so if we define the
matrix $M$ by
$$M =\begin{bmatrix}
 \theta_1 & \theta_2 & \theta_3 \\
\theta_1^{\ell} & \theta_2^{\ell} & \theta_3^{\ell} \\
 \theta_1^{m} & \theta_2^{m} & \theta_3^{m}
\end{bmatrix},
$$
then the equation $Mx=0$ has a nontrivial solution, and so $\det M =0$. Because $\theta_1\theta_2\theta_3 \neq 0$, this
is equivalent to
$(\theta_2^{\ell-1}-\theta_1^{\ell-1})(\theta_3^{m-1}-\theta_1^{m-1})=(\theta_3^{\ell-1}-\theta_1^{\ell-1})(\theta_2^{m-1}-\theta_1^{m-1}).$
Now we define the function $f$ by $f(x)=\frac{\theta_1^{m-1}-x^{m-1}}{\theta_1^{\ell-1}-x^{\ell-1}}$, so then
$f(|\theta_2|)=f(|\theta_3|)$ (note that $\ell$ and $m$ are odd). We claim however that $f$ is strictly increasing on
the interval $(0,\theta_1)$, which will give the required contradiction. To prove the claim, consider the function $g$
defined by $g(x)=f'(x)(\theta_1^{\ell-1}-x^{\ell-1})^2/x^{\ell-2}$. Then it is straightforward to show that
$g(x)=(m-\ell)x^{m-1}-(m-1)\theta_1^{\ell-1}x^{m-\ell}+(\ell-1)\theta_1^{m-1}$ and
$g'(x)=(m-1)(m-\ell)(x^{m-2}-\theta_1^{\ell-1}x^{m-\ell-1})$. So $g(\theta_1)=0$ and $g'(x)<0$ if $0<x<\theta_1$.
Therefore $g(x)>0$ if $0<x<\theta_1$, and this implies indeed that $f$ is strictly increasing on $(0,\theta_1)$. \epf

We tried to find (relevant) nonzero solutions to \eqref{eq:lswr} for some small $\ell>3$. We could not find integer
solutions, but we did find some solutions for $\ell=5$ with two conjugate noninteger eigenvalues and one integer
eigenvalue (for example $5,\frac12(-3+\sqrt{281}),\frac12(-3-\sqrt{281})$). We however have no clue whether there are
any regular four-eigenvalue graphs with these particular eigenvalues (note that there is still freedom in ``choosing''
the number of vertices and the valency). It thus remains open whether there are strongly $\ell$-walk-regular graphs
with $\ell>3$ that are not strongly $m$-walk-regular for all odd $m \neq \ell$.

We note that given $\ell$, and two of the three eigenvalues, say $\theta_2$ and
$\theta_3$, one can consider \eqref{eq:lswr} as a polynomial equation in
$\theta_1$, and this has two trivial roots $\theta_2$ and $\theta_3$. By using
again Lemma \ref{lem:descartes}, it follows that there can be no other real
root for $\ell$ even (which confirms part of Theorem \ref{thm:3or4}), and one
other real root for $\ell$ odd. In other words, for a given odd $\ell$, and
given $\theta_2$ and $\theta_3$, there is one real solution to \eqref{eq:lswr}
for $\theta_1$ that is different from $\theta_2$ and $\theta_3$. But typically
this solution seems to be unfit as an eigenvalue of a graph.

\section{The exceptional case}\label{sec:mu0}

What remains in this paper is to classify the strongly $\ell$-walk-regular
graphs with parameters $(\lal,0,\nul)$, that is, those with $\mul=0$. In this
case, we have the equation $A^{\ell}-\lal A- \nul I=O$ (by Lemma \ref{lem:A}).
Thus, $\G$ is strongly $\ell$-walk-regular graph with parameters
$(\lal,0,\nul)$ if and only if every eigenvalue of $\G$ is a root of the
polynomial $x^{\ell}-\lal x-\nul$. We again use Lemma \ref{lem:descartes} to
obtain the following classification.

\begin{proposition}\label{prop:mu=0}
Let $\ell>1$ and let $\G$ be a nonempty strongly $\ell$-walk-regular graph with
parameters $(\lal, 0, \nul)$. Then either $\G$ is a disjoint union of complete
graphs of the same order, or $\ell$ is odd and $\G$ is a disjoint union of
complete bipartite graphs of the same size and isolated vertices.
\end{proposition}

\pf First, let $\ell$ be even. Then by Lemma \ref{lem:descartes}, the
polynomial $x^{\ell}-\lal x-\nul$ has at most two real roots, so $\G$ has at
most two distinct eigenvalues. This implies that $\G$ is the disjoint union of
complete graphs of the same order (number of vertices) \cite[Thm. 6.4]{cds82}.

Next, let $\ell$ be odd. In this case, it follows that $\G$ has at most three distinct eigenvalues. Let us assume that
$\G$ is not a disjoint union of complete graphs of the same order, then $\G$ has three distinct eigenvalues. Consider
first the case that $\nul=0$. Then the eigenvalues of $\G$ are $\theta$, $0$, and $-\theta$, for some $\theta$. This
implies that $\G$ is indeed a disjoint union of complete bipartite graphs of the same size (number of edges) and
possibly some isolated vertices (cf. \cite[Thm. 6.5]{cds82}).

Finally, we consider the case where $\G$ has three distinct eigenvalues and $\nul>0$ (and $\ell$ is odd). Now $0$ is
not an eigenvalue, and we will show that this gives a contradiction. Also, in this case $\G$ cannot have two complete
graphs of different orders as components (because these would give different values for $\nul$). We assumed that $\G$
is not a disjoint union of complete graphs of the same order, so there is a component $C$ of $\G$ that is not a
complete graph. This component also has at most three eigenvalues, and so it has diameter at most $2$. Consider
vertices $u$ and $v$ in $C$ that are not adjacent. Then $u$ and $v$ are at distance two, and so there is a walk of
length $3$, and hence of length $\ell$, between $v$ and each of the neighbors of $u$. Because $\mul=0$, this implies
that $v$ is adjacent to all neighbors of $u$. Thus, every nonneighbor of $u$ (including $u$ itself) is adjacent to
every neighbor of $u$. It also follows easily that there are no edges among the neighbors of $u$ and that there are no
edges among the nonneighbors of $u$ (otherwise there would be a walk of length $3$ between $u$ and a nonneighbor of
$u$). Therefore $C$ is a complete bipartite graph, which gives the required contradiction. \epf

Note that the disjoint union of complete graphs $K_n$ is strongly regular, so
indeed it is strongly $\ell$-walk-regular for every $\ell$ (by Proposition
\ref{prop:srg}). If $\G$ is a disjoint union of complete bipartite graphs, all
of which have the same number of edges, and possibly some isolated vertices,
then indeed $\G$ is strongly $\ell$-walk-regular for every odd $\ell$. Note
finally that the disjoint union of complete graphs $K_2$ and $K_1$ (of
different orders) is a special case of a disjoint union of complete bipartite
graphs and isolated vertices, which is why we split the cases $\nul=0$ and
$\nul>0$ in the above proof.

\end{document}